\documentclass[11pt, a4paper]{article}
\usepackage[hmargin=2.5cm, vmargin=2cm, includefoot, twoside]{geometry}
\usepackage[backref=page,colorlinks,bookmarksopen=true]{hyperref}
\usepackage[latin1]{inputenc}  
\usepackage[T1]{fontenc}       
\usepackage[english]{babel}
\usepackage{amsmath}
\usepackage{amssymb}
\usepackage{amsthm}
\usepackage{latexsym}
\usepackage{mathrsfs}
\usepackage{graphics, graphicx}
\usepackage{tabularx}
\usepackage{multirow}
\usepackage{longtable}
\theoremstyle{plain}
\newtheorem{prop}{Proposition}

\newtheorem{lem}[prop]{Lemma}

\newtheorem{claim}{Claim}

\newtheorem*{thmMain}{Theorem}
\theoremstyle{definition}

\newtheorem*{remark}{Remark}
\theoremstyle{remark}

\newcommand{\C}{\mathbb{C}}

\newcommand{\R}{\mathbb{R}}
\newcommand{\Z}{\mathbb{Z}}

\newcommand{\la}{\langle}
\newcommand{\ra}{\rangle}
\newcommand{\inv}{^{-1}}

\newcommand{\g}{\mathfrak{g}}
\newcommand{\Iaff}{I_{\mathrm{aff}}}
\newcommand{\reDelta}{\; {}^{\mathrm{re}}\Delta}
\newcommand{\CAT}[1]{\mathsf{CAT}(#1)}%
\def\bs#1.{
              \def\temp{#1}
              \ifx\temp\empty
                   \mathcal{B}
              \else
                   \mathcal{B}(#1)
              \fi
}

\def\min{\mathop{\mathrm{min}}\nolimits}

\long\def\symbolfootnote[#1]#2{\begingroup%
\def\thefootnote{\fnsymbol{footnote}}\footnote[#1]{#2}\endgroup}
\begin{document}
\title%
{A uniform bound on the nilpotency degree of certain subalgebras of Kac-Moody algebras}
\author{Pierre-Emmanuel Caprace\thanks{F.N.R.S. research fellow, partially supported by the Wiener-Anspach Foundation.}}
\date{17th December 2006}

\maketitle

\begin{abstract}
Let $\mathfrak{g}$ be a Kac-Moody algebra and $\mathfrak{b}_1, \mathfrak{b}_2$ be Borel subalgebras of opposite
signs. The intersection $\mathfrak{b} = \mathfrak{b}_1 \cap \mathfrak{b}_2$ is a finite-dimensional solvable
subalgebra of $\mathfrak{g}$. We show that the nilpotency degree of $[\mathfrak{b}, \mathfrak{b}]$ is bounded
from above by a constant depending only on $\mathfrak{g}$. This confirms a conjecture of Y.~Billig and
A.~Pianzola \cite{BilligPia95}.
\end{abstract}

\newcommand{\subjclass}[1]{\symbolfootnote[0]{\noindent AMS classification numbers (2000):~#1.}}
\newcommand{\keywords}[1]{\symbolfootnote[0]{\emph{Keywords}:~#1}}


%
%
%
\section{Introduction}


Let $A$ be a generalized Cartan matrix, let $\mathfrak{g}$ be a Kac-Moody algebra of type $A$ with Cartan
decomposition $\mathfrak{g} = \mathfrak{h} \oplus \bigoplus_{\alpha \in \Delta} \mathfrak{g}_\alpha$ and $W$ be
its Weyl group. For each $w \in W$ we set $\Delta(w) = \{\alpha \in \Delta_+ \; | \; w.\alpha < 0\}$ and
$\mathfrak{g}_w = \la \mathfrak g_\alpha \; | \; \alpha \in \Delta(w) \ra$. It is known that $\Delta(w)$ is
finite and that $\mathfrak{g}_w$ is a finite-dimensional nilpotent subalgebra of $\mathfrak{g}$. The main result
of this paper is the following:

\begin{thmMain}\label{thm:nilp}
The nilpotency degree of $\mathfrak{g}_w$ is bounded from above by a constant depending on $A$, but not on $w$.
\end{thmMain}

This statement was conjectured in \cite[Conjecture~1]{BilligPia95}. In view of \cite[Proposition~1]{Ti87}, it is
equivalent to its group theoretic counterpart, which can be stated as follows. Let $G$ be the complex simply
connected Kac-Moody group of type $A$. For each $\alpha \in \reDelta$, let $U_\alpha$ be the one-parameter
subgroup of $G$ with Lie algebra $\g_\alpha$ and for all $w \in W$, let $U_w = \la U_\alpha \; | \alpha \in
\Delta(w) \ra$.

\begin{thmMain}
The nilpotency degree of $U_w$ is bounded from above by a constant depending on $A$, but not on $w$.
\end{thmMain}

The latter statement holds not only for the complex Kac-Moody group $G$, but in fact for any split or almost
split Kac-Moody group over an arbitrary field: this follows from a reformulation of the main theorem in terms of
root systems (see Proposition~\ref{prop:RootSyst} below) together with the description of commutation relations
in Kac-Moody groups \cite[Theorem~2]{Morita87}.

It is useful to keep in mind the group-side of the theory. For example, in the case where $A$ is of finite or
affine type, the preceding statement is an immediate consequence of the fact that the group $G$ is linear
(modulo center). On the other hand, for any other type of generalized Cartan matrix, the group $G$ is known to
be non-linear \cite[Theorem~7.1]{Cap05c}. Our proof of Theorem~\ref{thm:nilp} is based on a reduction to the
affine case. The main tools are, one the one hand, the classification of pairs of real roots whose sum is real
root, due to Y.~Billig and A.~Pianzola \cite{BilligPia95}, and on the other hand, on some sufficient conditions
on a set of roots to generate an affine subsystem, which were established in \cite{Ca06}.

\section{Nilpotent sequences in root systems}

\subsection{Definition}
A set of roots $\Phi \subset \Delta$ is called \textbf{prenilpotent} if there exist $w, w' \in W$ such that $w .
\Phi \subset \Delta_+$ and $w'.\Phi \subset \Delta_-$. In particular $\Phi \subset \; {}^{\mathrm{re}}\Delta$.
The set $\Phi$ is called \textbf{closed} if  for all $\alpha, \beta \in \Phi$ such that $\alpha + \beta$ is a
root, we have $\alpha + \beta \in \Phi$. Since the intersection of any collection of closed subsets is closed,
it makes sense to consider the closure of a set of roots. The following lemma is obvious:

\begin{lem}\label{lem:closure}
The closure of any prenilpotent set of roots is prenilpotent.
\end{lem}
\begin{proof}
The closure of a set $\Phi$ is contained in $(\sum_{\alpha \in \Phi} \Z_+ \alpha) \cap \Delta$.
\end{proof}

A sequence of roots $(\beta_k)_{k=1, \dots, n}$ is called \textbf{nilpotent} if it satisfies the following
conditions:
\begin{description}
\item[(NS1)] The set $\{\beta_1, \dots, \beta_k\}$ is prenilpotent.

\item[(NS2)] For each $k = 1, \dots, n$, we have $\sum_{i=1}^k \beta_i \in \Delta$.
\end{description}

For all $\alpha, \beta \in \Delta$, we have $[\g_\alpha, \g_\beta] \subset \g_{\alpha + \beta}$. Furthermore for
each $w \in W$, the prenilpotent set $\Delta(w)$ is closed. Therefore, the subalgebra $\g_w$ splits as a direct
sum $\g_w = \bigoplus_{\alpha \in \Delta(w)} \mathfrak{g}_\alpha$. It is easy to deduce from these basic facts
that the nilpotency degree of $\g_w$ coincides with the maximal possible length of a nilpotent sequence of roots
contained in $\Delta(w)$. Therefore, the following statement is equivalent to the main theorem and can be viewed
as its root system version:

\begin{prop}\label{prop:RootSyst}
The supremum of the set of lengths of nilpotent sequences of roots in $\Delta$ is finite.
\end{prop}

The proof of Proposition~\ref{prop:RootSyst} is deferred to Sect.~\ref{sect:proof}. We first need to collect a
series of subsidiary results: this is the purpose of the next subsection.

\subsection{On infinite root systems and their geometric realizations}

We freely use the standard notation and terminology on infinite root systems which can be found in
\cite[Chapter~5]{MP95}. We view $\Delta$ as the root system of a set of root data $\mathscr{D} = (A, \Pi,
\Pi^\vee, V, V^\vee, \la \cdot, \cdot\ra )$ over $\R$ and denote by $W$ its Weyl group. We assume that the
generalized Cartan matrix $A$ is finite. Moreover, we need to consider a geometric realization of $\reDelta$; we
henceforth denote by $X$ the interior of the Tits cone $\mathfrak X^\vee \subset V^\vee$. Recall that $X$ is
$W$-invariant and that the induced action is properly discontinuous \cite[Ch.~5, Prop.~15]{MP95}. For each root
$\alpha \in \reDelta$, we set $D(\alpha) = \{x \in X \; | \;  \la \alpha, x \ra > 0\}$ and $\partial \alpha =
\{x \in X \; | \; \la \alpha, x \ra = 0\}$. The set $\partial \alpha$ is called a \textbf{wall}; it is the trace
on $X$ of a hyperplane of $V^\vee$ and it cuts $X$ into two nonempty convex open cones, called
\textbf{half-spaces}, namely $D(\alpha)$ and $D(-\alpha)$. Note that walls and half-spaces are convex. The
notion of convexity will be crucial to our purposes.

\begin{remark}
Instead of the interior of the Tits cone, we might equally use the Davis complex associated with the Weyl group
$W$. This also provides a convenient geometric realization of $\reDelta$, which has no linear structure but is
instead equipped with a $W$-invariant $\CAT 0$-metric. This allows to define walls and half-spaces and yields an
appropriate notion of (geodesic) convexity. The $1$-skeleton of the Davis complex is nothing but the Cayley
graph of $W$ with respect to its canonical Coxeter generating set $S$. This graph may be embedded in the
interior of the Tits cone by considering as vertex set the $W$-orbit of a point in the interior of the
fundamental Weyl chamber and this makes it easy to pass from one viewpoint to the other. In the present note, we
keep the Tits cone viewpoint throughout, but we will be led to quote references which use rather the Davis
complex as a preferred geometric realization.
\end{remark}

The following lemma collects a few basic facts on pairs of roots:
\begin{lem}\label{lem:pairs}
Let $\alpha, \beta \in \reDelta$.
\begin{itemize}
\item[(i)] The subsystem generated by $\alpha$ and $\beta$ is finite if and only if $\partial \alpha$ meets
$\partial \beta$.

\item[(ii)]  $D(\alpha) \subset D(\beta)$ or $D(\alpha) \supset D(\beta)$ if and only if $\la \alpha, \beta^\vee
\ra \la \beta, \alpha^\vee \ra \geq 4 $ and $\la \alpha, \beta^\vee \ra >0$.

\item[(iii)] The pair $\{\alpha, \beta\}$ is prenilpotent if and only if $D(\alpha) \subset D(\beta)$, or
$D(\alpha) \supset D(\beta)$, or $\alpha \neq \pm \beta$ and the subsystem generated by $\alpha$ and $\beta$ is
finite.
\end{itemize}
\end{lem}
\begin{proof}
(i) follows from \cite[Ch.~5, Prop.~14]{MP95} and the fact that any finite subgroup of $W$ fixes a point of $X$.

\smallskip (ii). We may assume $\alpha \neq \pm \beta$, otherwise the desired assertion is obvious. In that
case, we have $D(\alpha) \subset D(\beta)$ or $D(\alpha) \supset D(\beta)$ only if the subsystem generated by
$\alpha$ and $\beta$ is infinite in view of (i). Since $\alpha \neq \pm \beta$, this in turn is equivalent to
$\la \alpha, \beta^\vee \ra \la \beta, \alpha^\vee \ra \geq 4 $. Now, if $\la \alpha, \beta^\vee \ra <0$, then
$\la \beta, \alpha^\vee \ra <0$ by \cite[Ch.~5, Prop.~8]{MP95} and it readily follows that $D(\alpha) \cap
D(\beta) \subset D(r_\beta(\alpha)) \cap D(r_\alpha(\beta))$. On the other hand, if $D(\alpha) \subset
D(\beta)$, then, transforming by $r_\beta$, we obtain $D(r_\beta(\alpha)) \subset D(-\beta)$ whence
$D(r_\beta(\alpha)) \cap D(\alpha) = \varnothing$. Similarly, if $D(\alpha) \supset D(\beta)$ then
$D(r_\alpha(\beta)) \cap D(\beta)=\varnothing$. This shows that if $D(\alpha) \subset D(\beta)$ or $D(\alpha)
\supset D(\beta)$ then $\la \alpha, \beta^\vee \ra >0$. The converse statement follows because, if $D(\alpha)
\not \subset D(\beta)$ and $D(\alpha) \not \supset D(\beta)$, then $D(-\alpha) \subset D(\beta)$ or $D(-\alpha)
\supset D(\beta)$.

\smallskip (iii). We may assume that the subsystem generated by $\{\alpha, \beta\}$ is infinite, otherwise the desired
assertion is easy. For any root $\alpha \in \reDelta$, we have $\alpha > 0$ if and only if the half-space
$D(\alpha)$ contains the Weyl chamber. Now the claim readily follows.
\end{proof}

\begin{lem}\label{lem:BilligPia}
There exists a constant $K$, depending only on the generalized Cartan matrix $A$, such that the following
condition holds. Given a prenilpotent pair $\{\alpha, \beta\} \subset \reDelta$ such that $\alpha + \beta$ is a
root, then $\la \alpha, \beta^\vee \ra \leq K$.
\end{lem}
\begin{proof}
Follows from \cite[Proposition~1 and Theorem~1]{BilligPia95}.
\end{proof}

\begin{lem}\label{lem:NR}
For any integer $n$, there exists a constant $L(n)$, depending on the generalized Cartan matrix $A$, such that
any prenilpotent set of at least $L(n)$ roots contains a subset $\{\alpha_1, \dots, \alpha_n\}$ of cardinality
$n$ such that $D(\alpha_1) \subsetneq D(\alpha_2) \subsetneq \dots \subsetneq D(\alpha_n)$.
\end{lem}
\begin{proof}
It is shown in \cite[Lemma 3]{NR03} that there exists a constant $L(2)$ such that any set of more than $L(2)$
walls contains a pair of parallel walls (i.e. non-intersecting walls). Combining this with Ramsey's theorem, it
follows that for any integer $n$, there exists a constant $L(n)$ such that any set of more than $L(n)$ walls
contains a set of $n$ pairwise parallel walls. Let now $\Phi$ be a prenilpotent set of roots of cardinality
greater than $L(n)$. Hence $\Phi$ contains a subset $\Phi_0$ of cardinality $n$ such that the elements of
$\partial \Phi_0 = \{\partial \alpha \; | \; \alpha \in \Phi_0\}$ are pairwise parallel. Since $\Phi_0$ is
prenilpotent, it follows from Lemma~\ref{lem:pairs}(iii) that the elements of $\{D(\alpha) \; | \; \alpha \in
\Phi_0 \}$ are totally ordered by inclusion. Thus they form a chain, as desired.
\end{proof}

\begin{lem}\label{lem:triangles}
There exists a constant $M$, depending on the generalized Cartan matrix $A$, such that the following property
holds. Let $\alpha, \alpha', \beta_0, \dots, \beta_n \in \reDelta$  be real roots such that:
\begin{itemize}
\item[(1)] The subsystem generated by $\{\alpha, \alpha', \phi_0\}$ is finite of rank~$2$. (Equivalently: we
have $\varnothing \neq \partial \alpha \cap \partial \alpha' \subset \partial \beta_0$.)

\item[(2)] $D(\beta_0) \subsetneq D(\beta_1) \subsetneq \dots \subsetneq D(\beta_n)$.

\item[(3)] For each $i=1, \dots, n$, the subsystem generated by $\{\alpha, \beta_i\}$ (resp. $\{\alpha',
\beta_i\}$) is finite. (Equivalently: the wall $\partial \beta_i$ meets both $\partial \alpha$ and $\partial
\alpha'$.)
\end{itemize}
If $n \geq M$, then the subsystem generated by $\{\alpha, \alpha', \beta_0, \dots, \beta_n\}$ is of irreducible
affine type; furthermore, it is contained in (a conjugate of) a parabolic subsystem of affine type of $\Delta$.
In particular the subsystem generated by $\{\beta_0, \dots, \beta_n\}$ is of  affine type and rank~$2$.
\end{lem}
\begin{proof}
See \cite[Theorem~8]{Ca06}.
\end{proof}

We will need to appeal to Lemma~\ref{lem:triangles} several times in the proof of
Proposition~\ref{prop:RootSyst}. The following lemma will be helpful when checking that the hypotheses of
Lemma~\ref{lem:triangles} are satisfied.

\begin{lem}\label{lem:Prenilp4}
We have the following:
\begin{itemize}
\item[(i)] Let $\{\alpha, \alpha', \gamma\} \subset \Delta$ be a prenilpotent set such that $D(\alpha)
\subsetneq D(\alpha')$ or $D(\alpha) \supsetneq D(\alpha')$ and $\la \alpha, \gamma^\vee \ra <0$. If $\partial
r_\gamma(\alpha)$ meets $\partial \alpha'$, then so does $\partial \gamma$.

\item[(ii)] Let $\{\alpha, \alpha', \beta, \beta'\} \subset \Delta$ be a prenilpotent set such that $D(\alpha)
\subsetneq D(\alpha')$ and $D(\beta) \subsetneq D(\beta')$. If $\partial \beta$ and $\partial \beta'$ both meet
$\partial \alpha'$ and if $\partial \beta'$ meets $\partial \alpha$, then $\partial \beta$ meets $\partial
\alpha$. Similarly, if $\partial \beta$ and $\partial \beta'$ both meet $\partial \alpha$ and if $\partial
\beta$ meets $\partial \alpha'$, then $\partial \beta'$ meets $\partial \alpha'$.
\end{itemize}
\end{lem}
\begin{proof}
(i). Up to replacing $\{\alpha, \alpha', \gamma\}$ by the prenilpotent set $\{-\alpha, -\alpha', -\gamma\}$, we
may assume without loss of generality that $D(\alpha') \subsetneq D(\alpha)$. Since $\la \alpha, \gamma^\vee \ra
<0$, it readily follows that $D(\alpha) \cap D(\gamma) \subset D(r_\gamma(\alpha))$. Assume now that $\partial
\gamma$ does not meet $\partial \alpha'$. Then we have $D(\alpha') \subsetneq D(\gamma)$ by
Lemma~\ref{lem:pairs}(iii) because, in view of Lemma~\ref{lem:pairs}, the wall $\partial \gamma$ meets $\partial
\alpha$ and moreover we have  $D(\alpha') \subsetneq D(\alpha)$. It follows that $D(\alpha') \subset D(\alpha)
\cap D(\gamma) \subsetneq D(r_\gamma(\alpha))$. In particular the wall $\partial \alpha'$ does not meet
$\partial r_\gamma (\alpha)$.

\smallskip
(ii). Suppose that $\partial \beta$ and $\partial \beta'$ both meet $\partial \alpha'$ and that $\partial
\beta'$ meets $\partial \alpha$. Now assume in order to obtain a contradiction that $\partial \beta$ does not
meet $\partial \alpha$. Then, by Lemma~\ref{lem:pairs} we have $D(\alpha) \subset D(\beta)$ or $D(\alpha)
\supset D(\beta)$. Since $\partial \beta$ meets $\partial \alpha'$ and $D(\alpha) \subset D(\alpha')$, we have
in fact $D(\alpha) \subset D(\beta)$. On the other hand, since $\partial \beta'$ meets $\partial \alpha$, it
follows that $D(\alpha) \cap D(-\beta')$ is nonempty. We deduce $\varnothing \neq D(\alpha) \cap D(-\beta')
\subset D(\alpha) \subset D(\beta)$. This contradicts the fact that $D(\beta) \subset D(\beta')$ which implies
$D(\beta) \cap D(-\beta') = \varnothing$.

The other assertion follows by considering the prenilpotent set $\{-\alpha, -\alpha', -\beta, -\beta'\}$.
\end{proof}

The next lemma provides useful sufficient conditions which ensure that a root belongs to a given parabolic
subsystem of affine type:

\begin{lem}\label{lem:affine}
Let $\Phi \subset \Delta$ be a parabolic subsystem of affine type and let $\alpha \in \Delta$. Then we have
$\alpha \in \Phi$ provided that one of the following conditions is fulfilled:
\begin{itemize}
\item[(i)] There exist $\beta, \beta' \in \Phi$ such that $D(\beta) \subsetneq D(\alpha) \subsetneq D(\beta')$.

\item[(ii)] There exist $\beta_1, \dots, \beta_8, \gamma \in \Phi$ such that the elements of $\{\partial
\beta_1, \dots, \partial \beta_8\}$ are pairwise parallel, the wall $\partial \alpha$ meets $\partial \beta_i$
for all $i = 1, \dots, 8$ and $\la \alpha, \gamma^\vee \ra \neq 0$.

\item[(iii)] There exist $\beta_1, \dots, \beta_n \in \Phi$ such that $D(\alpha) \subsetneq D(\beta_1)
\subsetneq \dots \subsetneq D(\beta_n)$ and $\la \beta_n, \alpha^\vee \ra < \frac{n}{2}$.
\end{itemize}
\end{lem}
\begin{proof}
(i) follows from \cite[Proposition~16 and Lemma~17]{Ca06}.

\smallskip (ii) follows from \cite[Lemma~11, Proposition~16 and Lemma~22]{Ca06}.

\smallskip (iii). Assume that $\beta_1, \dots, \beta_n \in \Phi$ are roots such that $D(\alpha) \subsetneq D(\beta_1)
\subsetneq \dots \subsetneq D(\beta_n)$ and suppose that $\alpha \not \in \Phi$. We must prove that $\la
\beta_n, \alpha^\vee \ra \geq \frac{n}{2}$.

Since $\Phi$ is of affine type, the condition $D(\beta_1) \subsetneq \dots \subsetneq D(\beta_n)$ implies that
the group $\la r_{\beta_1}, \dots, r_{\beta_n} \ra$ is infinite dihedral \cite[Theorem~D and
Proposition~14]{Ca06} and, hence, the subsystem $\Phi_0$ generated by $\{\beta_1, \dots, \beta_n\}$ is affine of
rank~$2$. Let $\beta_0 \in \Phi_0$ be the root such that $\{-\beta_0, \beta_1\}$ is a basis of $\Phi_0$, and
define inductively $\beta_{-k} =- r_{\beta_{-k+1}}(\beta_{-k+2})$ for all $k > 0$. Thus, for all $k \geq 0$, we
have $\beta_{-k} \in \Phi$, $D(\beta_{-k}) \subset D(\beta_{-k+1})$ and $\{-\beta_{-k}, \beta_{-k+1}\}$ is a
basis of $\Phi_0$. Note moreover that for all $k > 0$, we have $\la  \beta_k, \alpha^\vee \ra > 0$ by
Lemma~\ref{lem:pairs}(ii). We claim that there exists $k \geq 0$ such that $\la \beta_{-k}, \alpha^\vee \ra \leq
0$.

Suppose the contrary. If the wall $\partial \alpha$ meets at least $8$~elements of $\{\partial \beta_{-k} \; |
\; k> 0\}$, then we obtain $\alpha \in \Phi$ by (ii). Therefore, we may assume that $\partial \alpha$ is
parallel to almost every element of $\{\partial \beta_{-k} \; | \; k> 0\}$. Suppose now that $\partial
\beta_{-k} \subset D(\alpha)$ for some $k >0$. Since $\la  \beta_{-k}, \alpha^\vee \ra > 0$, we have
$D(\beta_{-k}) \subset D(\alpha) \subset D(\beta_1)$ by Lemma~\ref{lem:pairs}(ii) and, hence, we obtain $\alpha
\in \Phi$ by (i). Therefore, we may assume that $D(\alpha) \subset D(\beta_{-k}) \subset D(\beta_1)$ for almost
every $k >0$. Let $x \in \partial \alpha$ and $y \in \partial \beta_1$ be any points. It follows from the above
that the segment $[x, y]$ meets $\partial \beta_{-k}$ for almost every $k \geq 0$, which contradicts
\cite[Ch.~5, Prop.~6 and~7]{MP95}. This proves the claim.

Let now $i = \min \{ k \geq 0 \; | \;  \la \beta_{-k},\alpha^\vee \ra \leq 0\}$. Hence $\{-\beta_{-i},
\beta_{-i+1} \}$ is a basis of $\Phi_0$ and moreover $\la \beta_{-i}, \alpha^\vee \ra \geq 0$ and  $\la
\beta_{-i+1}, \alpha^\vee \ra \geq 1$. Define $\phi_1 = \beta_{-i+1}$, $\phi'_1 = r_{\beta_{-i+1}}(\beta_i)$ and
for all $k > 0$, set $\phi_{k+1} = r_{\beta_{-i+1}}r_{\beta_{-i}}( \phi_k)$ and $\phi'_{k+1}=
r_{\beta_{-i+1}}r_{\beta_{-i}}( \phi'_k)$. Thus we have $D(\beta_{-i+1}) = D(\phi_1) \subsetneq D(\phi'_1)
\subsetneq D(\phi_2) \subsetneq D(\phi'_2) \subsetneq \dots$ and moreover $\{\phi \in \Phi_0 \; | \;
D(\beta_{-i+1}) \subset D(\phi)\} = \{\phi_k, \phi'_k \; | \; k > 0\}$. In particular, we have $\beta_n \in
\{\phi_k, \phi'_k\}$ for some $k \geq \frac{n}{2}$.

Since $\{-\beta_{-i}, \beta_{-i+1} \}$ is a basis of $\Phi_0$, we may write $\phi_k = x_k.(-\beta_{-i}) +
y_k.\beta_{-i+1}$ for some non-negative integers $x_k, y_k$; similarly $\phi'_k = x'_k.(-\beta_{-i}) +
y'_k.\beta_{-i+1}$ for some $x'_k, y'_k \in \Z_+$. Since $\Phi_0$ is of affine type, an easy computation shows
that $y_{k+1}-y_k$ (resp. $y'_{k+1}-y'_k$) is a constant positive integer (it is independent of $k$). In
particular, the sequence $(y_k)$ (resp. $(y'_k)$) is a linear function of $k$ with positive integral coefficient
and, hence, we have $y_k \geq k$ (resp. $y'_k \geq k$) for all $k > 0$. It follows that
$$\begin{array}{rcl}
\la \phi_k, \alpha^\vee \ra & = & x_k.\la -\beta_{-i}, \alpha^\vee \ra + y_k. \la \beta_{-i+1}, \alpha^\vee
\ra\\
& \geq & y_k\\
& \geq & k.
\end{array}$$
Similarly, we obtain $\la \phi'_k, \alpha^\vee \ra \geq k$. Since $\beta_n \in \{\phi_k, \phi'_k\}$ for some $k
\geq \frac{n}{2}$, we finally obtain $\la \beta_n, \alpha^\vee \ra \geq \frac{n}{2}$
\end{proof}

\subsection{Proof of Proposition~\ref{prop:RootSyst}}\label{sect:proof}

Assume, in order to obtain a contradiction, that $\Delta$ contains nilpotent sequences of arbitrarily large
length. Therefore, by Lemma~\ref{lem:NR}, given any integer $n$, there exists a nilpotent sequence
$(\beta_i)_{i=1, \dots, k}$ such that the set $\{D(\beta_1), \dots, D(\beta_k)\}$ contains a chain of
half-spaces of length~$n$. Assume now that $k$ is fixed and let $I \subset \{1, \dots, k\}$ be a subset \emph{of
maximal possible cardinality} such that  $\{D(\beta_i) \; | \; i \in I\}$ is a chain of half-spaces. Provided
that $k$ is large enough, we may assume that $|I|$ is arbitrarily large.

For every $i \in \{1, \dots, k\}$, we set
$$I(i) = \{j \in I \; | \; j < i\}$$
and
$$J(i) = \{j \; | \; 1 \leq  j < i,\ \la \beta_j, \beta_i^\vee \ra < 0\}.$$

Note that for all $j \in J(i)$, the wall $\partial \beta_j$ meets $\partial \beta_i$ by Lemma~\ref{lem:pairs},
since the pair $\{\beta_i, \beta_j\}$ is prenilpotent.

\begin{claim}\label{cl:1}
For each $i \in I$, we have $|J(i)| \geq \frac{|I(i)|-K}{3}$, where $K$ is the constant of
Lemma~\ref{lem:BilligPia}.
\end{claim}

By (NS1) and Lemma~\ref{lem:BilligPia}, we have $ \sum_{j=1}^{i-1} \la\beta_j, \beta_i^\vee \ra \leq K$. On the
other hand, for each $j \in I$ we have $\la \beta_{j}, \beta_i^\vee \ra \geq 1$ by Lemma~\ref{lem:pairs}(ii) and
for all $j \in J(i)$, we have $\la \beta_j, \beta_i^\vee \ra \geq -3$. We deduce successively:
$$\begin{array}{rcl}
K & \geq &  \sum_{j \in I(i)} \la \beta_j, \beta_i^\vee \ra + \sum_{j < I, j \not \in I} \la \beta_j, \beta_i^\vee \ra \\
& \geq &  |I(i)| + \sum_{j \in J(i)} \la \beta_j, \beta_i^\vee \ra \\
& \geq & |I(i)| - 3|J(i)|.
\end{array}
$$\qed

\begin{claim}\label{cl:3}
Let $m$ be an integer such that $|I| > 4m$. Suppose that $I$ possesses an element $i$ such that $|J(i)| \geq
L(4m)$ and  $m > M$, where $L$ (resp. $M$) is as in Lemma~\ref{lem:NR} (resp. Lemma~\ref{lem:triangles}). Then
there is a set $\Iaff \subset \{1, \dots, k\}$ of cardinality $m$, such that $\{D(\beta_i) \; | \; i \in
\Iaff\}$ is a chain of half-spaces and the subsystem generated by $\{\beta_i \; | \; i \in \Iaff\}$ is of affine
type and rank~$2$, and is contained in a parabolic subsystem of affine type.
\end{claim}

By assumption, there exist $4m$ elements $\lambda_1, \dots, \lambda_{4m} \in J(i)$ such that
$D(\beta_{\lambda_1}) \subsetneq \dots \subsetneq D(\beta_{\lambda_{4m}})$. Let ${\gamma_1} = \beta_i + \la
\beta_i, \beta_{\lambda_1}^\vee \ra \beta_{\lambda_1}$ and ${\gamma_{4m}} = \beta_i + \la \beta_i,
\beta_{\lambda_{4m}}^\vee \ra \beta_{\lambda_{4m}}$. Note that, by Lemma~\ref{lem:closure}, the set $\{\beta_1,
\dots, \beta_k\} \cup \{\gamma_1, \gamma_{4m}\}$ is prenilpotent. Let also $I_- = \{ j \in I \; | \; D(\beta_j)
\subsetneq D(\beta_i)\}$ and $I_+ = \{ j \in I \; | \; D(\beta_i) \subsetneq D(\beta_j)\}$.
%

\smallskip
Suppose first that there exists a subset $I_0 \subset I_-$ of cardinality $m$ such that for each $j \in I_0$,
the wall $\partial \beta_j$ meets $\partial {\gamma_1}$. Since $\lambda_1 \in J(i)$ and since the pair
$\{\beta_i, \beta_{\lambda_1}\}$ is prenilpotent, it follows from Lemma~\ref{lem:pairs} that $\{\beta_i,
\beta_{\lambda_1}\}$ generates a finite subsystem. Furthermore, by Lemma~\ref{lem:Prenilp4}(i), the wall
$\partial \beta_{\lambda_1}$ meets $\partial \beta_j$ for all $j \in I_0$. Therefore, Lemma~\ref{lem:triangles}
ensures that $\{\beta_j \; | \; j \in I_0\}$ generates an affine subsystem of rank~$2$ which is contained in a
parabolic subsystem of affine type. Thus we are done in this case.

\smallskip
Suppose similarly that there exists a subset $I_0 \subset I_+$ of cardinality $m$ such that for each $j \in
I_0$, the wall $\partial \beta_j$ meets $\partial {\gamma_{4m}}$. Then, by the same argument as in the preceding
paragraph, we conclude that $\{\beta_j \; | \; j \in I_0\}$ generates an affine subsystem of rank~$2$ which is
contained in a parabolic subsystem of affine type.  Thus we are done in this case as well.

\smallskip
Suppose now that there exists a subset $I_1 \subset I_-$ of cardinality $m$ such that for each $j \in I_1$ and
for each $j' \in \{2, \dots, m\}$, the wall $\partial \beta_j$ meets $\partial \beta_{\lambda_{j'}}$. If for
some $j' \in \{2, \dots, m\}$, the wall $\partial \gamma_1$ is parallel to $\partial \beta_{\lambda_{j'}}$, then
it follows from Lemma~\ref{lem:Prenilp4}(ii) applied to $\{\beta_i, \beta_j, \beta_{\lambda_1}, \gamma_1\}$ that
$\partial \gamma_1$ meets $\partial \beta_j$ for all $j \in I_1$. Thus we are reduced to a case which has
already been settled. Therefore, we may assume that $\partial \gamma_1$ meets $\partial \beta_{\lambda_{j'}}$
for all $j'\in \{2, \dots, m\}$. In that case, Lemma~\ref{lem:triangles} implies that $\{\beta_{\lambda_{1}},
\dots, \beta_{\lambda_{m}}\}$ generates a subsystem of affine type and rank~$2$ which is contained in a
parabolic subsystem of affine type.

\smallskip
Suppose similarly that there exists a subset $I_1 \subset I_+$ of cardinality $m$ such that for each $j \in I_1$
and for each $j' \in \{3m+1, \dots, 4m-1\}$, the wall $\partial \beta_j$ meets $\partial \beta_{\lambda_{j'}}$.
Then, by the same argument as in the preceding paragraph using $\gamma_{4m}$ instead of $\gamma_1$, we conclude
that $\{\beta_{\lambda_{3m+1}}, \dots, \beta_{\lambda_{4m}}\}$ generates an affine subsystem of rank~$2$ which
is contained in a parabolic subsystem of affine type.   Thus we are done in this case as well.

\smallskip
Let us now define $I' \subset I$ to be the subset consisting of all those $j$'s such that $\partial \beta_j$
meets $\partial \beta_{\lambda_{j'}}$ for some $j' \in \{m+1, \dots, 3m\}$. By Lemma~\ref{lem:Prenilp4}(ii), if
$j \in I' \cap I_-$ then $\partial \beta_j$ meets $\partial \beta_{\lambda_{j'}}$ for all $j' \in \{1, \dots,
m\}$. It follows that we may assume $|I' \cap I_- | < m$, otherwise we are reduced to a case which has already
been settled. Similarly, if $j \in I' \cap I_+$ then $\partial \beta_j$ meets $\partial \beta_{\lambda_{j'}}$
for all $j' \in \{3m+1, \dots, 4m\}$ and, as above, we may assume that $|I' \cap I_+ | < m$. Since $I = I_- \cup
\{i\} \cup I_+$, it follows that $I'= (I' \cap I_-) \cup \{i\} \cup (I' \cap I_+)$ and, hence, the last case
which remains to be treated is when $|I'| < 2m$. Note that, by definition, the set $\{\partial \beta_j \; | \; j
\in I \backslash I'\} \cup \{\partial \beta_{\lambda_{m+1}}, \dots,
\partial \beta_{\lambda_{3m}}\}$ consists of pairwise parallel walls. Therefore, the set $\{D(\beta_j) \; | \; j
\in I \backslash I'\} \cup \{D(\beta_{\lambda_{m+1}}), \dots, D(\beta_{\lambda_{3m}})\}$ is a chain of
half-spaces of length $|I| - |I'| + 2m > |I|$. This contradicts the maximality property of $I$, thereby showing
that this last case does not occur.\qed

\begin{claim}
Suppose that there exists a set $\Iaff = \{\lambda_1, \dots, \lambda_n\} \subset \{1, \dots, k\}$ of cardinality
$n > 6.L(8)+2K + 6$ such that $D(\beta_{\lambda_1}) \subsetneq \dots \subsetneq D(\beta_{\lambda_n})$ and
$\{\beta_{\lambda_1}, \dots, \beta_{\lambda_n}\}$ generates an affine subsystem of rank~$2$ which is contained
in a parabolic subsystem of affine type, where $L$ is the function of Lemma~\ref{lem:NR}. Then there exists a
nilpotent sequence $(\beta'_j)_{j=1, \dots, k'}$, such that the set $\{\beta'_1, \dots, \beta'_{k'}\}$ contains
$\{\beta_{\lambda_j}\; | \; j = x+1, \,  x+2, \, \dots, \, n-x\}$, where $x = 6.L(8)+2K +6$, and is contained in
a parabolic subsystem of affine type of $\Delta$.
\end{claim}

We make the following definitions:
$$i := \min \{\lambda_j \; | \; j = x+1, \, x+2, \dots, \, n-x\}, \hspace{1cm} k':= k-i+2,$$
$$\beta'_1 := \beta_i, \hspace{1cm} \beta'_2 := \sum_{j=1}^{i-1} \beta_j  \hspace{1cm}
\text{and}   \hspace{1cm}  \beta'_j := \beta_{i+j-2} \text{ for all } j = 3, \dots, k'.$$ The set $\{\beta'_1,
\dots, \beta'_{k'}\}$ is prenilpotent by Lemma~\ref{lem:closure}, thus the sequence $(\beta'_j)_{j \leq k'}$
satisfies (NS1). It also satisfies (NS2) because so does $(\beta_j)_{j \leq k}$. Hence $(\beta'_j)_{j \leq k'}$
is a nilpotent sequence. Furthermore, it follows from the definition that the set $\{\beta'_1, \dots,
\beta'_{k'}\}$ contains $\{\beta_{\lambda_j}\; | \; j = x+1, \, x+2, \, \dots, \, n-x\}$. Let $\Phi \subset
\Delta$ be the parabolic subsystem of affine type containing $\{\beta_{\lambda_1}, \dots, \beta_{\lambda_n}\}$.
We now show by induction on $m$ that $\{\beta'_j \; |j=1, \dots, m \}$ is contained in $\Phi$. This is true by
hypothesis for $m=1$, hence the induction can start.

Let $m>1$. Note that by induction, we have $\sum_{j=1}^{m-1} \beta'_j \in \Phi$ and moreover $\la
\sum_{j=1}^{m-1} \beta'_j , {\beta'_m}^\vee \ra \neq 0$ by (NS2).  Therefore, if $\partial \beta'_m$ meets at
least $8$~elements of $\{\partial \beta_{\lambda_1}, \dots,
\partial \beta_{\lambda_n}\}$, then we have $\beta'_m \in \Phi$ by Lemma~\ref{lem:affine}(ii) and we are done.
We henceforth assume that $\partial \beta'_m$ meets at most $7$~elements of $\{\partial \beta_{\lambda_1},
\dots,
\partial \beta_{\lambda_n}\}$.

If $\partial \beta'_m$ meets some element of $\{\partial \beta_{\lambda_8}, \dots,
\partial \beta_{\lambda_{n-7}}\}$, then by the above the triple $\{\partial \beta_{\lambda_1}, \partial
\beta'_m, \partial \beta_{\lambda_n}\}$ consists of pairwise parallel walls since walls are convex. Moreover,
since $\partial \beta_{\lambda_j} \subset D(\beta_{ \lambda_n}) \backslash D(\beta_{\lambda_1})$ for all $j = 2,
\dots, n-1$ and since $\partial \beta'_m$ meets $\partial \beta_{\lambda_j}$ for some such $j$, we obtain
$\partial \beta'_{m} \subset D(\beta_{ \lambda_n}) \backslash D(\beta_{\lambda_1})$. Finally, we deduce from
Lemma~\ref{lem:pairs}(iii) that $D(\lambda_1) \subsetneq D(\beta'_m) \subsetneq D(\lambda_n)$. By
Lemma~\ref{lem:affine}(i), this implies that $\beta'_m \in \Phi$.

It remains to consider the case when $\partial \beta'_m$ meets no element of $\{\partial \beta_{\lambda_8},
\dots, \partial \beta_{\lambda_{n-7}}\}$. Thus the set $\{D(\beta'_m), D(\beta_{\lambda_8}), \dots,
D(\beta_{\lambda_{n-7}})\}$ is a chain of half-spaces by Lemma~\ref{lem:pairs}(iii) and we may assume that
either $D(\beta'_m) \subsetneq D(\beta_{\lambda_8})$ or $D(\beta'_m) \supsetneq D(\beta_{\lambda_{n-7}})$
otherwise we may conclude again using Lemma~\ref{lem:affine}(i). Define $J'(m) := \{ j < m \; | \; \la \beta'_j
, {\beta'_m}^\vee \ra <0\}$. Note that $1 \not \in J'(m)$ otherwise $\partial \beta'_m$ would meet $\partial
\beta'_1$ by Lemma~\ref{lem:pairs}. Note that for all $j \in J'(m)$, we have $\la \beta'_j, {\beta'_m}^\vee \ra
\geq - 3$ in view of Lemma~\ref{lem:pairs}. Therefore, since $\la \sum_{j=1}^{m-1} \beta'_j , {\beta'_m}^\vee
\ra \leq K$ by (NS2) and Lemma~\ref{lem:BilligPia}, we deduce:
$$\begin{array}{rcl}
\la \beta'_1, {\beta'_m}^\vee \ra & \leq & K - \sum_{j=2}^{m-1} \la \beta'_j, {\beta'_m}^\vee \ra\\
& \leq & K - \sum_{j \in J'(m)} \la \beta'_j, {\beta'_m}^\vee \ra \\
& \leq & K + 3 \; |J'(m)|.
\end{array}$$
If $K + 3 \, |J'(m) | < \frac{x-6}{2}$, then we obtain $\beta'_m \in \Phi$ by Lemma~\ref{lem:affine}(iii), as
desired. Otherwise, we have $|J'(m)| \geq \frac{x-2K-6}{6} = L(8)$ by the definition of $x$. Therefore, the set
$J'(m)$ contains $8$~elements $j_1, \dots, j_8$ such that the walls $\partial \beta'_{j_1}, \dots,
\partial \beta'_{j_8}$ are pairwise parallel. By Lemma~\ref{lem:pairs}, the wall $\partial \beta'_m$ meets
$\partial \beta'_j$ for each $j \in J'(m)$, since the pair $\{\partial \beta'_j, \partial \beta'_m\}$ is
prenilpotent. Furthermore, we have $\{ \beta'_j \; | \; j \in J'(m)\} \subset \Phi$ by induction. We finally
conclude that $\beta'_m \in \Phi$ by Lemma~\ref{lem:affine}(ii).\qed

\medskip

We are now ready to obtain a final contradiction. The above claims show that the existence of nilpotent
sequences of arbitrarily large length in $\Delta$ implies the existence of nilpotent sequences of arbitrarily
large length, entirely contained in parabolic subsystems of affine type of $\Delta$. Note that there are only
finitely many orbits of such subsystems under the Weyl group action. Thus there must exist nilpotent sequences
of arbitrarily large length, entirely contained in some fixed parabolic subsystem of affine type of $\Delta$. As
mentioned in the introduction, this is impossible because it contradicts the fact that Kac-Moody groups of
affine type are linear modulo center. Here are some more details.

Linearity of affine Kac-Moody groups follows from their well known realization as matrix groups over rings of
Laurent polynomials: if $G$ is the complex simply connected Kac-Moody group of untwisted affine type $X^{(1)}$
(notation of \cite[Ch.~4]{Ka90}) and $\mathbf{G}$ denotes the simple simply connected algebraic group scheme of
type $X$, then there is a central homomorphism $\varphi : G \to \mathbf{G}(\C[t, t\inv])$ by
\cite[Sect.~7.3]{Ti81/82}. Using the divisibility of the subgroup $U_w < G$, it follows that the Zariski closure
of $\varphi(U_w)$ in $ \mathbf{G}(\overline{\C[t, t\inv]})$ is connected. Since it is moreover nilpotent, it is
contained in a Borel subgroup of $\mathbf{G}$. Therefore, the nilpotency degree of $U_w$ is bounded from above
by the solvability degree of Borel subgroups of $\mathbf{G}$, which is of course independent of $w$. A similar
argument also applies to the case of twisted affine groups, since these can be viewed as almost split forms of
untwisted affine groups \cite[Sect.~7.3]{Ti81/82} and, hence, are linear as well.

An alternative way to prove that the length of nilpotent sequences in a root system $\Delta$ of affine type is
uniformly bounded, is to use the description of real roots in $\Delta$ provided by \cite[Proposition~6.3]{Ka90}.
This result shows that if $(\beta_j)_{j=1, \dots, k}$ is a nilpotent sequence of $\Delta$, then $(\overline
\beta_j)_{j=1, \dots, k}$ is a nilpotent sequence of the finite (possibly non-reduced) root system
$\overline{\reDelta}$, where $\alpha \mapsto \overline \alpha$ denotes  the orthogonal projection introduced in
\cite[\S6.2]{Ka90}. In particular, every nilpotent sequence of $\Delta$ is of length at most
$|\overline{\reDelta}|$.\qed

\providecommand{\bysame}{\leavevmode\hbox to3em{\hrulefill}\thinspace}
\providecommand{\MR}{\relax\ifhmode\unskip\space\fi MR }
\providecommand{\MRhref}[2]{%
  \href{http://www.ams.org/mathscinet-getitem?mr=#1}{#2}
} \providecommand{\href}[2]{#2}

\noindent
D\'epartement de Math\'ematiques\\
Universit\'e Libre de Bruxelles, CP 216\\
Boulevard du Triomphe\\
1050 Bruxelles, Belgique.\\

\noindent
\emph{Current address:}\\
Mathematical Institute,\\
University of Oxford,\\
24--29 St Giles',\\
Oxford OX1 3LB,\\
United Kingdom.\\
{\tt caprace@maths.ox.ac.uk}

\end{document}